\input amstex
\documentstyle{amsppt}
\document
\magnification=1200
\NoBlackBoxes
\nologo
\pageheight{18cm}


\bigskip

\centerline{\bf THE NOTION OF DIMENSION}
\smallskip

\centerline{\bf IN GEOMETRY AND ALGEBRA\footnotemark1}
\footnotetext{Based on the talks delivered at the AMS annual meeting, Northwestern U., Oct. 2004,
and Blythe Lectures, U. of Toronto, Nov. 2004.}

\medskip

\centerline{\bf Yuri I. Manin}

\medskip

\centerline{\it Max--Planck--Institut f\"ur Mathematik, Bonn, Germany,}

\centerline{\it and Northwestern University, Evanston, USA}

\bigskip

{\bf Abstract.} This talk reviews some mathematical and physical
ideas related to the notion of dimension.
After a brief historical introduction, various modern constructions
from fractal geometry, noncommutative geometry, and theoretical physics
are invoked and compared.

\bigskip

\hfill{{\it Glenn Gould disapproved of his own recording of Goldberg variations.}}

\hfill{{\it ``There is a lot of piano playing going on there, and I}}

\hfill{{\it mean that as the most disparaging comment possible.''}}
\smallskip

\hfill{{\it NYRB, Oct. 7, 2004, p.10}}

\bigskip

\centerline{\bf \S 0. Introduction}

\medskip

{\bf 0.1. Some history.} The notion of dimension belongs to the most fundamental
mathematical ideas. In the Western civilization and school system,
we become pretty early exposed to the assertion that
the dimension of  our physical space is three
(and somewhat later, that time furnishes the fourth dimension).

\smallskip

However, what does such a statement actually mean?

\smallskip

The mental effort needed to grasp the meaning
of ``three'' in this context is qualitatively different
from the one involved in making sense of
a sentence like ``There are three chairs in this room''.
Counting dimensions, we are definitely {\it not}
counting ``things''.
 
\smallskip

Even if we make a great leap to abstraction, and accept 
Cantor's sophisticated definition of a whole number
as a cardinality of a finite set, life does not get much easier.
A Cantorian set is supposed to be 
``any collection of definite,
distinct objects  of our perception  or our thought.''
But what are exactly  these ``distinct objects'',
these ``non--things'' which we take from our physical space
and project into our mind?

\smallskip

Euclid (ca 300 BC), as some of great thinkers before and after him,
taught us not to bother so much about what things ``are''
but rather how to think about them orderly {\it and} creatively.

\smallskip

Let us reread his often quoted passages where dimension is indirectly
involved. I use the delightfully archaic rendering of [He]:

\medskip

{\it From BOOK I, On plane geometry:

\smallskip

1. A point is that which has no part.

2. A line is breadthless length.

3. The extremities of a line are points.

[...]

5. A surface is that which has length and breadth only.

6. The extremities of a surface are lines.

\medskip

From BOOK XI, On spatial geometry:

\smallskip

1. A solid is that which has length, breadth, and depth.

2. An extremity of a solid is a surface.}

\medskip

Two observations are immediate. First, Euclid directs
our imagination to ``kitchen physics'': e.g.
we are supposed to grasp right away what is ``a part''
and to have no difficulty to imagine an entity
without parts. He describes semantics of basic geometric notions
in terms of a slightly refined non--verbal everyday experience. 

\smallskip

Second, for a possible discoverer and a 
great practitioner of the axiomatic method,
he is strangely oblivious about some one--step
logical implications of his definitions. If we take the 
``extremity'' (I will also use the modern term boundary)
of a solid ball, it must be a surface, namely, sphere.
Now, the boundary of this surface is {\it not} a line,
contrary to Book I, 6, because
it is empty!
 
\smallskip

Thus, Euclid misses a great opportunity here:
if he stated the principle

\medskip

\centerline{\it ``The extremity of an
extremity is empty'',} 

\medskip

he could be considered
as the discoverer of the 

\medskip

\centerline{BASIC EQUATION OF HOMOLOGICAL ALGEBRA:} 
$$
d^2=0.
$$
For a historian of culture, the reason of this strange
blindness is obvious: ``emptiness'' and ``zero'' as legitimate
notions, solidly built into systematic scientific thinking,
appear much later. 

\smallskip

Even contemporary thought
periodically betrays persistent intellectual
uneasiness regarding emptiness: compare M.~Heidegger's
reification of ``Nothing'', and various versions  of ``vacuum state''
in the quantum field theory. 

\smallskip

Skipping two milleniums, turn now to Leibniz. The following
excerpt from [Mand], p. 405, introduces a modern algebraic
idea related to our further discussion of fractional dimensions:

\medskip

``[...] the idea of fractional integro--differentiation
[...] occurred to Leibniz, as soon as he has developed his version of calculus and invented
the notations $d^kF/dx^k$ and $(d/dx)^kF$. In free
translation of Leibniz's letter to de l'H\^opital dated September 30, 1695
[...]:

\smallskip 

{\it Johann Bernoulli seems to
have told you of my having mentioned to him
a marvelous analogy which makes it possible to say in a way that
successive differentials are in geometric progression.  One can ask what would be a differential
having as its exponent a fraction. [...] Although this seems removed
from Geometry, which does not yet know of such fractional
exponents, it appears that one day these paradoxes will yield useful consequences,
since there is hardly a paradox without utility.}''

\medskip

I read the first part of this quotation as a reference
to what became later known as Taylor series: in the formula
$f(x+dx)=\sum_{n=0}^{\infty} \dfrac{f^{(n)}(x)}{n!} (dx)^n$
the consecutive terms are ``successive differentials in geometric progression''.

\smallskip

By extension, I will interpret Leibniz's  quest
as follows: ``Make sense of the formal expression
$f(x) (dx)^s$'', with arbitrary rational, or real, or even
complex value of $s$ (we may add $p$--adic values as well).

\smallskip

Nowadays it is easy to give a Bourbaki--style answer to this quest.

\smallskip

 Let $M$ be a differentiable manifold,
$s$  an arbitrary complex number. 

\medskip

Then we can construct:

\medskip

(i) A rank one complex vector bundle $V_s$ {\it of $s$--densities}  on $M$,
which is trivialized over each coordinate neighborhood $(x_i)$
and for which the transition multiplier from $(y_j)$ to $(x_i)$  is 
$$
|\,\roman{det}\,(\partial x_i/\partial y_j)|^s
$$

\smallskip

(ii) Its sections locally can be written as $h(x) |dx|^s$.
Spaces of sections $W_s$ with various integrability,
differentiability, growth etc conditions are called {\it $s$--densities.}

\smallskip

(iii) (Some) sections of $V_1$ are measures, so that they can be
integrated on $M$.

\smallskip

(iv) This produces a scalar product on (various spaces of) densities:
$$
W_s\times W_{1-s} \to \bold{C}:\quad (f|dx|^s, g|dx|^{1-s})\mapsto
\int_M fg |dx|.
$$

\medskip

Roughly speaking, subsets $N\subset M$ of (normalized) fractional dimensions $s\,\roman{dim}\,V \in \bold{R}$
appear when we learn that $W_s$ can be {\it integrated} along them,
as $W_1$ can be integrated along $M$.
Various ramifications and amplifications of this idea will be
reviewed below. 

\medskip

Euclid and Leibniz were chosen to represent in this Introduction 
respectively
right- and left- brain modes of thinking, characteristic 
traits of which became well known after popularizations
of Roger Sperry's work on brain asymmetry (Nobel Prize 1981).
In mathematical thinking, they roughly correspond to the
dichotomies Geometry/Algebra, Vision/Formal Deduction etc.
As I have written elsewhere ([Ma4]): 

\smallskip

``A natural or acquired
predilection towards geometric or algebraic
thinking and respective mental objects is often expressed
in strong pronouncements, like Hermann Weyl's
exorcising ``the devil of abstract
algebra'' who allegedly struggles with
``the angel of geometry''  for the soul
of each mathematical theory. (One is reminded of an even more sweeping
truth: {\it ``L'enfer -- c'est les autres''}.) 

\smallskip

Actually, the most fascinating thing about Algebra and Geometry
is the way they struggle to help each other
to emerge from the chaos of non--being, from
those dark depths of subconscious where all roots of
intellectual creativity reside. What one ``sees''
geometrically must be conveyed to others in words and symbols.
If the resulting text can never be a perfect
vehicle for the private and personal vision, the vision itself
can never achieve  maturity without
being subject to the test of written speech. The latter
is, after all, the basis of the social 
existence of mathematics. 

\smallskip

A skillful use of the interpretative algebraic language possesses 
also a definite therapeutic 
quality. It allows one to fight the obsession which
often accompanies contemplation of enigmatic
Rorschach's blots of one's  private imagination.

\smallskip

When a significant new unit of meaning (technically, 
a mathematical definition
or a mathematical fact) emerges from such a struggle,
the mathematical community spends some time
elaborating all conceivable implications
of this discovery. (As an example, imagine the development
of the idea of a continuous function, or a Riemannian metric,
or a structure sheaf.)
Interiorized, these implications prepare new firm ground for
further flights of imagination, and more often than not
reveal the limitations of the initial formalization
of the geometric intuition. 
Gradually the discrepancy between
the limited scope of this unit of meaning 
and our newly educated and enhanced geometric vision
becomes glaring, and the cycle repeats itself.''

\medskip

This all--pervasive Left/Right dichotomy has also a
distinctive social dimension, which recently led to the juxtaposition
of ``Culture of the Word'' and ``Culture of the Image''.
The society we live in becomes more and more dominated
by mass media/computer generated images
(to which visual representations of ``fractals'', sets of 
non--integer dimension, marginally belong).
Paradoxically, this technologically driven evolution
away from ``logocentrism'', often associated
with modernity and progress, by relying heavily upon
right brain mental faculties, projects us directly into
dangerously archaic states of collective consciousness.

\medskip

{\bf 0.2. Plan of the paper.} The first section presents several
contexts in which one can define dimension transcending the core
intuition of ``number of independent degrees of freedom'': 
Hausdorff--Besicovich dimension, dimensional regularization,
Murray--von Neumann dimension. 

\smallskip

The second section introduces dimension in supergeometry
and discusses the question: {\it what is dimension of $Spec\,\bold{Z}$?}

\smallskip

The third section is devoted to the spectrum of dimensions
Leibniz style arising in the theory of modular forms.

\smallskip

Finally, in the fourth section I review some recent constructions
introducing fractional dimensions in homological algebra.
Although their source was Mirror Symmetry, I have chosen to
present a theorem due to Polishchuk which puts
this subject in the context of noncommutative geometry and
Real Multiplication program for quantum tori.
With some reluctance, I decided to omit
another fascinating development, which can also be related 
to Mirror Symmetry,
``motives of fractional weights'' (cf. [And]). 

\medskip

{\it Acknowledgements.} After my talk at the AMS meeting,
Ed Frenkel reminded me about my old report [Ma1].
Matilde Marcolli has read 
the first draft of this paper and suggested
to include the discussion of Connes' notion of {\it
dimension spectrum,} and Lapidus' related notion
of {\it complex dimensions}. She has also written detailed instructions
to that effect. The whole \S 2 owes its existence 
to them, and  provides an additional
motivation for my presentation of the Lewis--Zagier theory in \S 3.
Sasha Beilinson and Sasha Polishchuk drew  my attention
to the old R.~MacPherson's explanation of perverse sheaves
on simplicial complexes. I am grateful to all of them.

\bigskip

\centerline{\bf \S 1. Fractional dimensions: a concise collector's guide}

\medskip

{\bf 1.1. A left/right balanced notion I: Hausdorff-Besicovich
dimension.} The scene here is a metric space $M$.
What is defined and counted: {\it HB--dimension} of an arbitrary subset $S\subset M$ with a compact closure.

\smallskip

Strategy of count:

\smallskip

(i) Establish that a Euclidean $d$--dimensional ball $B_{\rho}$ of radius $\rho$
has volume 
$$
vol_d(B_{\rho})= \frac{\Gamma(1/2)^d}{\Gamma(1+d/2)}\,\rho^d.
\eqno(1.1)
$$
Here $d$ is a natural number.

\smallskip

(ii) {\it Declare} that a $d$--dimensional ball $B_{\rho}$ of radius $\rho$
has volume given by the same formula
$$
vol_d(B_{\rho})= \frac{\Gamma(1/2)^d}{\Gamma(1+d/2)}\,\rho^d.
$$
for {\it any} real $d$.

\smallskip

(iii) Cover $S$ by a finite number of balls of radii $\rho_m$.

\smallskip

(iv) Make a tentative count: measure $S$ {\it as if}
it were $d$--dimensional for some $d$:
$$
v_d(S):=\roman{lim}_{\rho\to 0}\roman{inf}_{\rho_m<\rho} \,\sum_m vol_d
(B_{\rho_m}).
$$

\smallskip

(v) Grand Finale: There exists $D$ such that $v_d(S) =0$ for
$d>D$ and $v_d(S) =\infty$ for $d<D.$

\smallskip

{\it This $D$ is declared to be the HB--dimension of $S$.}  

\medskip

A set of non--integral Hausdorff--Besicovich dimension is 
called {\it a fractal} (B.~Mandelbrot).

\medskip

The general existence of $D$ is a remarkable mathematical phenomenon,
akin to those that appear in the description of phase transitions
and critical exponents in physics. To get a feeling
how it works, consider first the simple example: how do we
see that $[0,1]$ isometrically embedded in $M$ has dimension one?
Basically, we can cover $[0,1]$ by $N$ closed balls of diameter $\rho=N^{-1}$
centered at points of $S$.
In a tentative count assuming dimension $d$, we get
approximate volume $c_dN\cdot (1/2N)^d$ which tends to 0 
(resp. to $\infty$) with $N\to\infty$ when $d>1$ (resp. $d<1$.)
Hence $D=1$.

\smallskip

A similar counting (which can be easily remade into a formal proof)
shows that the classical Cantor subset $C\subset [0,1]$
has dimension $\roman{log}\,2/\roman{log}\,3$, so is a fractal.

\smallskip

A further remark: the constant involving gamma--factors in the
formula (1.1) for $vol_d(B_{\rho})$ does not influence the value of $D$.
However, some $S$ of HB--dimension $D$ may have a definite value of 
$v_D(S)$, that is, be HB--measurable. The value of this
volume will then depend on the normalization.
Moreover, one can slightly change the scene, and work with,
say, differentiable ambient manifolds $M$ and $s$--densities 
in place of volumes of balls. This leads to the picture of integration
of densities I referred to in the Introduction.

\smallskip

I do not know, whether some $S$ are so topologically ``good'' as to deserve the name of
$D$--dimensional manifolds. Is there a (co)homology theory
with geometric flavor
involving such fractional dimensional sets?

\medskip

{\bf 1.2. Fractional dimensions in search of a space: 
dimensional regularization of path integrals.}
Here is a very brief background; for details,
see [Kr], [CoKr] and references therein.

\smallskip

Correlators in a quantum field theory are given heuristically
by Feynmann path integrals.
A perturbative approach to defining such an integral
produces a formal series whose terms are indexed by Feynman
graphs and are familiar finite dimensional integrals.

\smallskip

However, each term of such a formal series usually
diverges. A procedure of regularizing it by
subtracting appropriate infinites is called  
{\it regularization/renormaliza\-tion.} Each such procedure 
involves a choice of a certain parameter, a value of a
mass scale, of an interaction constant, etc.
which is then made variable in such a way, that
integrals become finite at ``non--physical'' values of this
parameter. A study of their analytic behavior near the physical point
then furnishes concrete divergent terms, or {\it counterterms},
which are subtracted.

\smallskip

{\it Dimensional regularization} is a specific regularization procedure
which replaces the physical dimension of space time
(4 in the case of a scalar nonstringy field theory) by
a complex variable $D$ varying in a small neighborhood of 4.

\smallskip

Instead of extrapolating volumes of balls to non--integral $D$,
one extrapolates here the values of a Gaussian integral:
$$
\int e^{-\lambda} |k|^2 d^Dk=\left(\frac{\pi}{\lambda}\right)^{D/2}.
\eqno(1.2)
$$
However, no explicit sets making geometric sense of (1.2)
occur in the theory.

\smallskip

The germ of the $D$--plane near $D=4$ then becomes the base 
of a flat connection with an irregular singularity.
The regularization procedure can be identified with
taking the regular part of a Birkhoff decomposition
(A. Connes -- D. Kreimer). 

\smallskip

If one wishes to think of spaces of such such complex dimensions,
one should probably turn to noncommutative geometry:
cf. subsections 1.4 1nd 2.4 below.

\medskip

{\bf 1.3. A left/right balanced notion II: Murray--von Neumann factors.}
Here the introductory scene unfolds in 
a linear space $M$ (say, over $\bold{C}$.)
What is counted: dimension of a linear subspace $L\subset M$.

\smallskip

How fractional dimensions occur: if the usual linear
dimensions of $L$ and $M$ are both infinite, it might happen nevertheless that
a relative dimension $\dfrac{\roman{dim}\,L}{\roman{dim}\,M}$
makes sense and is finite.

\smallskip

To be more precise, we must first rewrite the finite
dimensional theory stressing the matrix algebra
$E_M:=\roman{End}\,M$ in place of $M$ itself.

\smallskip

Replace $L$ by (the conjugacy class of) projector(s) $p_L\in E_M$:
$p_L^2=p_L$, $p_L(M)=L$.

\smallskip

Construct a (normalized) trace functional $tr:\, E_M\to \bold{C}$.
A natural normalization condition here is $tr (id_M)=1$.
Another normalization condition might be $tr (p_L)=1$
where $L$ is a subspace having no proper subspaces.

\smallskip

Define the (fractional) dimension of $L$ as $tr\, (p_L)$.

\medskip

Murray and von Neumann pass from this elementary picture
to the following setting: take for $L$ a Hilbert space;
define $W$--algebras acting on $L$ as algebras of bounded operators
with abstract properties similar to that of $E_M$. 

\smallskip

For a reasonable class of such algebras (factors), study 
{\it the spectrum of fractional dimensions:}
values of (possibly normalized) trace functionals $t$ on the 
equivalence classes of projections $t (p_L)$
(projections are selfadjoint projectors).

\smallskip

The remarkably beautiful Murray and von Neumann classification
theorem then says that the normalized spectra of 
dimensions can be exactly of five types:

\smallskip

I${}_n$:\ $\{1,\dots ,n\}$.

\smallskip

I${}_{\infty}$:\ $\{1,\dots ,\infty \}$.

\smallskip

II${}_1$:\ $[0,1]$

\smallskip

II${}_{\infty}$:\ $[0,+\infty ]$

\smallskip

III:\ $\{0,+\infty \}$.

\medskip

{\bf 1.4. Moving to the left with Alain Connes: Noncommutative Geometry.}
A vast project (actually, a vast building site) of
Noncommutative Geometry is dominated by two different
motivations. A powerful stimulus is furnished by physics:
quantum theory replaces commuting observables by generally
noncommuting operators. 
\smallskip

Another motivation comes from
mathematics, and its conception  is definitely a
left brain enterprise. 
It is a fact that the study of all more or less
rigid geometric structures (excluding perhaps homotopical topology)
is based on a notion of (local) functions on the respective
space, which in turn evolved from the idea of coordinates
that revolutionized mathematics. After Grothendieck,
this idea acquired such a universality that, for example,
{\it any} commutative ring $A$ now comes equipped with a space
on which $A$ is realized as an algebra of functions:
the scheme $\roman{Spec}\,A$.

\smallskip

In Noncommutative Geometry we allow ourselves to think about
noncommutative rings (but also about much more general structures,
eventually categories, polycategories etc) as coordinate
rings of a ``space'' (resp. sheaves etc on this space). What constitutes the ``existential characteristics''
of such a space, what algebraic constructions reveal
its geometry, how to think about them orderly and creatively, --
these are the challenges that fascinate many practitioners
in this field.

\smallskip

One device for efficient training of our geometric
intuition is the consistent study of commutative geometry from noncommutative
perspective. In particular, an important role is played by certain
 spaces which appear as ``bad quotients''
of some perfectly sane commutative manifolds.  A standard example
is furnished by the space of leaves of a foliation: such a space is well defined
as a set but generally has a very bad topology. Alain Connes in [Co2]
sketched the general philosophy and provided a ghost of beautiful examples
of such situations. The starting point in many
cases is the following prescription for
constructing a noncommutative ring describing
a bad quotient $M/\Cal{R}$: take a function ring $A$
of $M$ and replace it by a certain crossed product 
$A\rtimes \Cal{R}.$

 \smallskip

Probably, the simplest example of a bad quotient is provided not
by a foliation, but by the trivial action of a, say, finite
group $G$ on a point. The quotient $S:=\{pt\}/G$ with respect
to such an action is represented by the group algebra
$A_S:= \bold{C}[G]$.  We imagine  $A_S$ as ``an algebra of functions on
a noncommutative space $S^{nc}$\,'' and apply to it the generic dictionary
of the {\it Algebra $\Leftrightarrow$ Geometry} correspondence:

\smallskip

a measure on $S^{nc}$ := a linear functional on $A_S$;

\smallskip

a vector bundle on $S^{nc}$ := a projective module over $A_S$ ...  etc. 

\medskip

Fractional dimensions of von Neumann and Murray
setting reappear in this context, with potentially very different
geometric interpretation. 

\smallskip 

At this point, it's worth stressing
the difference between classical fractals
and such noncommutative spaces : the former are embedded in ``good''
spaces, the latter are their projections. Living in a Platonic cave,
we have more psychological difficulties in recognizing and describing
these projections.

\medskip

{\bf 1.5. Digression on databases.} A large database $B$ with links can be imagined as a vast metric space. We can envision its graph approximation: pages $\Rightarrow$
vertices, links $\Rightarrow$ edges. Metric can be defined by the
condition that a link has length one;
or else: length of a link is the relative number of hits.

\smallskip

Approximate dimension $d$ then can be introduced: a (weighted) number of pages
accessible in $\le R$ links is approximately
$c R^{d}$.

\smallskip

Some experimental work with actual databases produces definitely
non--integer dimensions ([Ma]).

\smallskip

Databases are used for search of information. A search in $B$ usually produces a ``bad subset'' $S$
in $B$, for example, all contexts of a given word.

\smallskip

If the database $B$ registers results of scientific
observations (Human Genome project, cosmology), what we would like
to get from it, is instead (a fragment of) {\it a new scientific theory.}

\smallskip

Arguably, such a theory is rather {\it ``a bad quotient''} than a bad
subset of $B$.

\smallskip

Imagine all Darwin's observations registered as
a raw data base, and imagine how evolutionary theory might have  been
deduced from it: drawing bold analogies  and performing a drastic compression. Arguably, both
procedures are better modeled by ``bad equivalence relations''
than by contextual search. 

\bigskip

\centerline{\bf \S 2. Exotic dimensions}

\medskip

Up to now, I mostly used the word
``dimension'' in the sense ``the number of
degrees of freedom'' (appropriately counted).
In the title of this section, ``exotic
dimension'' means ``an unusual degree of freedom'',
and is applied, first, to odd dimensions of supergeometry
and second, to the arithmetical line $Spec\,\bold{Z}.$
Subsections 2.1 and 2.2 can be read as a post--scriptum to the review
written twenty years ago: see [Ma1] and [At].

\smallskip

Starting with 2.3, we explain the notion that
(doubled real parts of) zeroes and poles of various zeta--functions
of geometric origin $Z(X,s)$ can be viewed as a
``dimension spectrum''. One source of this notion is the
algebraic/arithmetic geometry of varieties over finite or number
fields. Another source, and the term ``dimension spectrum'',
is Connes' work on noncommutative Riemannian geometry (cf. [Co2], VI,
IV.3.$\gamma$)
and a closely related work of M.~Lapidus and collaborators
in fractal geometry ([LaPo], [LavF1], [LavF2]).

\medskip

{\bf 2.1. Supergeometry.} Coordinate rings in supergeometry
are $\bold{Z}_2$--graded and supercommutative: we have
$fg=(-1)^{|f||g|}gf$ where $|f|$ denotes the parity of $f$.
This Koszul sign rule applies generally in all algebraic
constructions. The algebraic coordinate ring of
an affine space $A^{m|n}_K$ over  a field $K$ is $K[x_1,\dots ,x_m;
\xi_1,\dots ,\xi_n]$ where $x$'s are even and $\xi$'s
are odd. Its (super)dimension is denoted $m|n$. 

\smallskip

Since odd functions are nilpotents, the usual intuition
tells us that, say, $A_K^{0|n}$ can be only imagined
as an ``infinitesimal neighborhood'' of the point
$Spec\,K$. This seemingly contradicts our desire to see this
superspace as a ``pure odd manifold''. A slightly more
sophisticated reasoning will convince us that
$A_K^{0|n}$ does have the defining property of a manifold:
its cotangent sheaf (universal target of odd derivations) is free,
again due to the Koszul sign rule.

\smallskip

Basics of all geometric theories can be readily extended to
superspaces. Deeper results also abound, in particular, Lie--Cartan classification
of simple Lie algebras is extended in a very interesting way.

\smallskip

The physical motivation for introducing odd coordinates
was Fermi statistics for elementary particles,
and the conjecture that laws of quantum field theory
include supersymmetry of appropriate Lagrangians.

\medskip

{\bf 2.2. What is the value of dimension of $Spec\,\bold{Z}$?}
 
{\it Answer 1: $\roman{dim}\,Spec\,\bold{Z}=1.$} This is the
common wisdom. Formally, one is the value of Krull
dimension of $\bold{Z}$, maximal length of a chain
of embedded prime ideals. Krull dimension can be viewed as a natural
algebraization of Euclid's inductive definition of dimension.
From this perspective, primes $p$ are zero--dimensional
points of $Spec\,\bold{Z}$, images of geometric points
$Spec\,\bold{F}_p\to Spec\,\bold{Z}.$

\smallskip

{\it Answer 2: $\roman{dim}\,Spec\,\bold{Z}=3.$}
One can argue, however, that $Spec\,\bold{F}_p$ is 
{\it not} zero--dimensional, because its fundamental
group, $Gal\,\overline{\bold{F}}_p/\bold{F}_p$,
is the same as (the completed) fundamental group
of the circle. Images of $Spec\,\bold{F}_p\to Spec\,\bold{Z}$
should be then visualized as {\it loops} in the space
$Spec\,\bold{Z}$, and one can define their linking numbers
which turn out to be related to the reciprocity laws and
Legendre and R\'edei symbols. For this reason,
$Spec\,\bold{Z}$ ``must'' be three--dimensional:
cf. a review in [Mor] and references therein.

\smallskip

More systematically, in \'etale topology
of  $Spec\,\bold{Z}$ one can observe 3--dimensional
Poincar\'e duality: see [Maz].

\smallskip

Before the advent of \'etale topology and even schemes,
Selberg made the remarkable discovery that
lengths of closed geodesics in hyperbolic spaces
behave like primes: Selberg's zeta functions are
close relatives of Riemann's zeta. This gives additional
weight to the idea that primes ``are'' loops.

\smallskip

{\it Answer 3: $\roman{dim}\,Spec\,\bold{Z}=\infty$} ??
This guess involves the conjectural existence
of a geometrical world defined over 
``an absolute point''  $Spec\,\bold{F}_1$ where
$\bold{F}_1$ is a mythical field with one element.
For some insights about this world, see 
[Ti], {Sm1], [Sm2], [KapSm], [Ma2], [Sou].

\smallskip

In particular, Soul\'e in [Sou] defined a category
of varieties over  $Spec\,\bold{F}_1$ which presumably
should be thought of as varieties of finite type.
This category {\it does not contain} $Spec\,\bold{Z}$.
In fact, objects $V$ of this category are defined
via properties of their purported base extensions
$V\times_{\bold{F}_1}\bold{Z}$, whereas $\bold{Z}\times_{\bold{F}_1}\bold{Z}$ 
remains tantalizingly elusive.

\smallskip

And if $Spec\,\bold{Z}$ is not a finite type
object, it can hardly  have a finite dimension.

\medskip

{\bf 2.2.1. Summary.}
The discussion so far can be interpreted as leading
to the following conclusion: not only the arithmetical
degree of freedom  $Spec\,\bold{Z}$ is exotic,
but the value of the respective dimension
is not just a real number or infinity, but a new entity
which deserves a special attention. We probe the arithmetical degree
of freedom by studying its interaction with
``geometric'' degrees of freedom, in particular, studying algebraic
varieties over finite and number fields and rings.

\smallskip

For this reason, the whole arithmetic geometry, in particular,
Arakelov's insights and their subsequent development, 
Deninger's program [De1]--[De4], relations with noncommutative geometry as in [Co5],[ConsMar], and Haran's
visions ([Ha]), will bear upon our future
enlightened decision about what $\roman{dim}\,Spec\,\bold{Z}$ 
actually is. Below we will briefly remind the
role of zeta functions from this perspective.

\medskip

{\bf 2.3. Zeta functions and weights in arithmetic geometry.}
Let $C$ be a smooth irreducible
projective algebraic curve defined over $\bold{F}_q$.
Its zeta function can be defined by a Dirichlet series
and an Euler product, in perfect analogy with Riemann zeta:
$$
Z(C,s)=\sum_a\frac{1}{N(a)^s}=
\prod_x\frac{1}{1-N(x)^{-s}}\,.
\eqno (2.1)
$$
Here $x$ runs over closed points of $C$, playing role of primes,
and $a$ runs over effective cycles rational over $\bold{F}_q$.
It is an elementary exercise to rewrite (2.1) in terms
of a generating function involving all numbers
$\roman{card}\,V(\bold{F}_{q^f}).$ The latter can
be interpreted as the numbers of fixed points of powers
of the Frobenius operator $Fr$ acting upon $C(\overline{\bold{F}}_q)$.
A. Weil's remarkable insight consisted in postulating
the existence of a cohomology theory and a Lefschetz type formula counting
these fixed points, proving it for curves, and
conjecturally extending it to general projective manifolds
over finite fields. The result for curves reads
$$
Z(C,s)=\prod_{w=0}^2 \roman{det}\left((\roman{Id} - Fr\cdot q^{-s})\,|\,
H^w(C)\right)^{(-1)^{w-1}} =
$$
$$
=\prod_{w=0}^2 Z(h^w(C),s)^{(-1)^{w-1}}\, .
\eqno(2.2)
$$
More generally, for a smooth irreducible projective
manifold $V$ defined over $\bold{F}_q$ we can define
the zeta function by a formula similar to (2.1) and,
after Grothendieck and Deligne, prove
the formula
$$
Z(V,s)=\prod_{w=0}^{2\roman{dim}\,V} \roman{det}\left((\roman{Id} - Fr\cdot q^{-s})\,|\,
H^w(V)\right)^{(-1)^{w-1}} =
$$
$$
=\prod_{w=0}^{2\roman{dim}\,V} Z(h^w(V),s)^{(-1)^{w-1}}\,.
\eqno(2.3)
$$
Here $H^w(V)$ denotes \'etale cohomology of weight $w$,
whereas $h^w(V)$ refers to the motivic piece of $V$ of weight $w$
which is a kind of universal cohomology.
According to the Riemann--Weil conjecture proved by Deligne,
the roots $\rho$ of $Z(h^w(V),s)$ lie on the vertical line
$\roman{Re}\,\rho =\dfrac{w}{2}.$

\smallskip

Thus we can read off the spectrum of dimensions in which
$V$ ``manifests itself'' nontrivially (i.e. by having a nontrivial
cohomology group) by

\smallskip

(i) Counting fixed points of Frobenius on $V(\overline{\bold{F}}_q)$.

\smallskip

(ii) Looking at the zeroes and poles of the zeta function
produced by this count.

\smallskip

Note also that the expression occurring in (2.3) very naturally
appears in the supergeometry of the total cohomology space
$H^*(V)$ graded by the parity of weight: this is simply
the inverse superdeterminant of the operator 
$\roman{Id} - Fr\cdot q^{-s}$. More generally,
quantum cohomology introduces a quite nontrivial and nonlinear
structure on this cohomology considered as a supermanifold
and not just $\bold{Z}_2$--graded linear space.
This is how supergeometry enters classical mathematics
from the back door.

\smallskip

Let us return now to $Spec\,\bold{Z}.$ Deninger in [De1]
suggested to write the following analog of (2.2)
for Riemann's zeta multiplied by the $\Gamma$--factor interpreted as
the Euler factor at arithmetic infinity:
$$
Z(\overline{Spec\,\bold{Z}},s):=2^{-1/2}\pi^{-s/2}\Gamma \left(\frac{s}{2}\right)
\zeta (s)= \frac{\prod_{\rho}\frac{s-\rho}{2\pi}}{\frac{s}{2\pi}\frac{s-1}{2\pi}} =
$$
$$
= (?)\,\prod_{\omega =0}^2
\roman{DET}\left(\frac{s\cdot\roman{Id}-\Phi}{2\pi}\,|\,H^{\omega}_?
(\overline{Spec\,\bold{Z}})\right)^{(-1)^{w-1}}\, .
\eqno(2.4)
$$
Here the notation $\prod_{\rho}$ as well as the conjectural $\roman{DET}$
refers to the ``zeta regularized'' infinite products which are defined by
$$
\prod_i \lambda_i:= \roman{exp}\,\left(-\frac{d}{dz}
\sum_i\lambda_i^{-z}\left|_{z=0}\right.\right)\,.
\eqno(2.5)
$$
The second equality sign in (2.4) is a theorem, whereas the
last equality sign expresses a conjecture about the existence
of some cohomology theory and a Frobenius type operator on it.
In fact, $\Phi$ should be considered rather as a logarithm of 
Frobenius: the direct comparison must be made between (2.4)
and (2.2) rewritten with the help of another identity
$$
1-\mu q^{-s}=\prod_{\alpha : q^{\alpha}=\mu} \frac{\roman{log}\,q}{2\pi i}
\,(s-\alpha ).
$$
Finally, to study  the interaction between $Spec\,\bold{Z}$
and geometric dimensions,
one considers zeta functions of schemes of finite type
over $Spec\,\bold{Z}$ which are, say, models of
smooth projective manifolds over number fields; or even motives
of this type. A series of partial results and sweeping conjectures
suggests a similar picture for such zetas, with
real parts of zeroes/poles producing a spectrum of
``absolute weights'' of arithmetical schemes.

\smallskip

For more details, see a discussion in [Ma3],
and more recent speculations on the nature of
geometry behind the Frobenius $\Phi$ and Deninger's cohomology.
In particular, [De3], [De4] postulate existence of dynamical
systems underlying this geometry, whereas [ConsMar]
introduces  noncommutative spaces responsible for 
$\Gamma$--factors for curves, and Connes in [Co5] outlines an approach
to Riemann hypothesis by way of noncommutative geometry.

\bigskip

{\bf 2.4. Dimension spectra of spectral triples.} Zeta functions of another
kind arise in the context of Connes {\it spectral triples}.
This is a reformulation of basic data of Riemannian geometry which can
be directly generalized to the noncommutative case: cf. [Co3], [CoMos].

\smallskip

We will now briefly describe the constructions relevant
to our discussion of dimension spectra.

\smallskip

A spectral triple consists of data $(A,H,D)$, where $A$ is a
$*$--subalgebra of bounded operators on a Hilbert space $H$,
whereas $D$ is an unbounded self--adjoint operator on $H$, with compact
resolvent. The relevant compatibility condition between $A$ and $D$
reads as follows: commutators $[a,D]$ are bounded for all
$a\in A$. 

\smallskip

The prototype of such a structure is the triple
$(C^\infty(X),L^2(S),D)$ associated to a compact spin manifold
$X$, where  $S$ is the spinor bundle and $D$ the Dirac operator.
Examples of spectral triples in genuinely noncommutative cases have
been constructed in the context of quantum groups
([Co6], [DaLSSV]), in arithmetic geometry ([ConsMar]), and
as a proposed geometric model of elementary particle physics
([Co4]).

\smallskip

Spectral triples of finite summability degree (that is, where $|D|^z$ is
trace class for some $z$) provide the stage for a universal local
index formula (Connes--Moscovici) for the cyclic cohomology Chern
character associated to the index problem
$$ Ind_D : K_*(A) \to \bold{Z}. $$
The local index formula is given in terms of the Wodzicki residue on
the algebra of pseudodifferential operators of $(A,H,D)$. Extending
the Wodzicki residue to this context amounts to extending the Dixmier
trace to operators of the form $b|D|^{-z}$, where $b$ is in
the algebra $B$ generated by the elements $\delta^n(a)$, for
$\delta(a)=[|D|,a]$ and $a\in A$. 

\smallskip

It is here that the notion of dimension spectrum naturally enters the
scene. Namely, the fact that the Wodzicki residue continues to make
sense and defines a trace depends upon the properties of a family of
zeta functions associated to the spectral triple $(A,H,D)$,
$$ 
\zeta_b (z)=\roman{Tr}\,(b\,|D|^{-z}), \ \ \  b\in B.  
$$
In particular, the dimension spectrum $\Sigma\subset \bold{C}$ is
the minimal subset such that, for all $b\in B$, the 
zeta function $\zeta_b$ extends holomorphically to $\bold{C}\smallsetminus
\Sigma$. 

\smallskip

If this subset is discrete and singularities
of $\zeta_b(s)$ are simple poles, 
one can extend the Dixmier trace by the formula 
$$ 
\int\, b  := \roman{Res}_{z=0}\, \roman{Tr}(b\,|D|^{-z}). 
$$
The case of the Connes--Moscovici local index formula where
multiplicities appear in the dimension spectrum is treated by applying
renormalization group techniques ([Co3], [CoMos]).

\smallskip

In what sense can one think of $\Sigma$ as a ``set of dimensions'' for the
noncommutative manifold $(A,H,D)$ ? 

\smallskip

In the context of this paper, a comparison with (2.2)--(2.4)
may help the reader. Connes includes in the dimension spectrum
only singularities of his zetas. It seems that they should
be compared with zeroes of (the numerator of) (2.2), (2.4),
that is, with poles of the inverted arithmetical zeta,
corresponding to the motivic part of {\it odd weight}.
This inversion corresponds in supergeometry to the
parity change, which might explain the remark made in the section
II.1, p. 205, of [CoMos], about advantages of treating the odd case.

\smallskip

If $(A,H,D)$ comes from a $p$--dimensional compact spin manifold,
the dimension spectrum of it is contained in $\{\,n\in\bold{Z}\,|\,
n\le p \}$ ([CoMos], p. 211), and the relevant singularities of zeta are simple poles.

\smallskip

One can associate spectral triples to certain fractal sets
and calculate their spectra. Even more straightforward constructions
can be given for special fractal sets like {\it fractal strings }
or {\it generalized fractal strings} of Lapidus
and van Frankenhuysen, cf. [LavF2], Ch. 3.

\smallskip

Sometimes it turns out that that the whole
spectrum lies on the line $\roman{Re}\, s=D_0$, where $D_0$
is the Hausdorff--Besicovich dimension. This again agrees
with the arithmetic geometry case: cf. our discussion
of Riemann hypothesis after formula (2.3). 
The extra factor $1/2$ appearing there is nicely explained by a peculiar normalization:
roughly speaking, in algebraic geometry dimension of complex line is 1, so that
dimension of the real line should be  $1/2$.

\smallskip

The dimension spectrum has the expected behavior with respect to the 
product of spectral triples. In arithmetic geometry,
similar arguments lead to the highly speculative
picture of Kurokawa's tensor product of zeta functions 
which presumably reflects a direct product operation
over ``the absolute point'' $\bold{F}_1$, cf. [Ma3].

\smallskip

The question of possible relations between dimension spectra, the
local index formula, and dimensional regularization arises naturally
in [CoMar2], where the Connes--Kreimer theory of perturbative
renormalization is reformulated as a Riemann--Hilbert problem for
a certain class of flat connections with irregular singularities
(equisingular). In fact, the universal singular frame, which
produces ``universal counterterms'' for all the renormalizable field
theories, has exactly the same rational coefficients that appear in
the local index formula of Connes--Moscovici. This suggests the 
possibility that the missing geometry underlying dimensional regularization
may be found in noncommutative spaces whose dimension spectrum lies
fully off the real line.

\smallskip

The fact that $Spec\,\bold{Z}$ ``manifests itself nontrivially'' both
in dimensions one and three, raises the possibility
that it can be connected with  the dimension spectrum of the (as
yet conjectural) spectral triple for a noncommutative
space such as the adelic quotient considered by Connes in the spectral
realization of the zeros of zeta ([Co5]), which is in turn the space of
commensurability classes of 1-dimensional $\bold{Q}$-lattices of [CoMar1].

\smallskip

In the next section I will review some recent results
due to J.~Lewis and D.~Zagier providing a geometric
scene for the interpretation of dimension spectrum of
a Selberg zeta function. Although spectral triples do not
explicitly appear there, the spirit is very similar.

\bigskip

\centerline{\bf \S 3. Modular forms and weights}

\medskip

{\bf 3.1. Classical modular forms.} A classical modular form
of weight $w+2$ with respect to a subgroup
$\Gamma\subset SL(2,\bold{Z})$ is a meromorphic
function $f$ on upper
half--plane $H$ satisfying
$$
f\left(\frac{az+b}{cz+d}\right)=f(z)(cz+d)^{w+2}
\eqno(3.1)
$$
for all fractional linear transformation from $\Gamma$.
Let us assume that $w$ is an even integer; then (3.1)
means that the formal expression  $f(z)(dz)^{(w+2)/2}$
is $\Gamma$--invariant and hence is
(the lift of) a Leibniz's higher differential on $X_{\Gamma}=\Gamma\setminus H$. 
\smallskip

One can also rewrite (3.1) differently, by 
looking at the universal elliptic curve $E_{\Gamma}\to X_{\Gamma}$
and the Kuga--Sato variety which is (a compactification of) $E_{\Gamma}^{(w)}:=E_{\Gamma}\times_{X_{\Gamma}}\dots\times_{X_{\Gamma}}
E_{\Gamma}$ ($w$ times). Namely, $E_{\Gamma}^{(w)}$
is the quotient of $\bold{C}^w\times H$ with respect to
to the group $\Gamma^{(w)}:=\bold{Z}^{2w}\rtimes \Gamma$
where the group of shifts $\bold{Z}^{2w}$ acts upon fibers by
$$
((t_k),z)\mapsto ((t_k+m_kz+n_k),z)
$$
whereas $\Gamma$ acts by
$$
((t_k),z)\mapsto \left(\left(\frac{t_k}{cz+d}\right),\frac{az+b}{cz+d}\right).
$$
In this notation, (3.1) means that the meromorphic
volume form $f(z) dz\wedge dt_1\wedge\dots\wedge dt_w$
comes from a volume form $F$ on $E_{\Gamma}^{(w)}$. Thus $w+1$ is the 
dimension of a classical space. Periods of cusp forms $f$, which are
by definition integrals $\int_0^{i\infty} f(z)z^kdz$, $0\le k\le w$,
can be also expressed as integrals of $F$ over appropriate cycles in 
$E_{\Gamma}^{(w)}$.

\smallskip

Below I will describe modular forms of {\it fractional weight}  of two types: 

\smallskip

(i) Serre's forms with $p$--adic weights;

\smallskip

(ii) Lewis--Zagier period functions with complex weights.

\medskip

{\bf 3.2. Serre's modular forms of $p$--adic weight.}
Here we will consider the case $\Gamma = SL(2,\bold{Z})$. Let
a modular form be finite at all cusps, in particular
given by its Fourier series of the form $\sum_{n=0}^{\infty} a_nq^n$,
$q:=e^{2\pi iz}.$ Former $w+2$ is now denoted $k$ and also called
a weight.

\medskip

{\bf 3.2.1. Definition.} {\it A $p$--adic modular form is a formal series
$$
f=\sum_{n=0}^{\infty} a_nq^n, \quad a_n\in \bold{Q}_p,
$$
such that there exists a sequence of modular forms $f_i$ 
of weights $k_i$ with rational coefficients $p$--adically
converging to $f$ (in the sense of uniform convergence of coefficients).}

\medskip

{\bf 3.2. Theorem.} {\it The $p$--adic limit
$$
k = k(f) := \roman{lim}_{i\to \infty} k_i \in \roman{lim}_{m\to \infty}
\bold{Z}/(p-1)p^m\bold{Z}=\bold{Z}_p\times \bold{Z}/(p-1)\bold{Z}
$$
exists and depends only of $f$.

\smallskip

It is called the $p$--adic weight of $f$.} 

\medskip

For a proof, see [Se]. This theorem enhances some constructions
which initially appeared in the ory of $p$--adic interpolation
of $L$--series; see also [Ka] for a broader context.

\smallskip 

Can one make sense of a $p$--adic limit of Kuga spaces $E_{\Gamma}^{(k_i-1)}$,
or rather, appropriate motives?  

\medskip

{\bf 3.3. Modular forms at the boundary.} Quotients 
$X_{\Gamma}=\Gamma\setminus H$ are noncompact modular curves.
In algebraic geometry, they are compactified by adding cusps $\Gamma\setminus \bold{P}^1(\bold{Q})$.
Recently in several papers it was suggested that
one should consider as well the ``invisible'' part of the modular
boundary consisting of the $\Gamma$--orbits of irrational
points in $\bold{P}^1(\bold{R}).$ The space 
$B_{\Gamma}:=\Gamma\setminus\bold{P}^1(\bold{R})$ is an archetypal bad quotient
which should be treated as a noncommutative space:
see [CoMar], [MaMar], [Ma4], and references therein.
Below we will discuss, what objects should be considered as
modular forms on $B_{\Gamma}.$

\smallskip

We will restrict ourselves by the basic case
$\Gamma =GL(2,\bold{Z}).$ There is another description of 
$B_{\Gamma}$: it is the set of the equivalence classes
of $\bold{R}$ modulo the equivalence relation
$$
x\equiv y \Leftrightarrow \exists\,m,n\ \roman{such\ that}\
T^mx=T^ny
\eqno(3.2)
$$
where the shift operator T is defined by
$$
T:\, x \mapsto \frac{1}{x}-\left[ \frac{1}{x} \right], 
$$
Consider instead the dual shift on functions given by the formal
operator $L=L_1$:
$$
(Lh)(x):=\sum_{k=1}^{\infty}
\frac{1}{(x+k)^{2}}\,
f\left(\frac{1}{x+k}
\right), 
\eqno(3.3)
$$
The meaning of dualization in this context is clarified
by the following formula:
$$ 
\int_{[0,1]} f \cdot Lh \, dx =
\int_{[0,1]} (f|T) \, h \, dx. 
$$
Generalizing (3.3), we can introduce the formal operator $L_s$ on functions $\bold{R}\to \bold{C}$:
$$
(L_sh)(x):=\sum_{k=1}^{\infty}
\frac{1}{(x+k)^{2s}}
h \left(\frac{1}{x+k}
 \right) .
$$
To put it more conveniently, consider
$$
L:=\sum_{k=1}^{\infty} \left[\matrix 0 & 1\\1 & k\endmatrix\right]\in \bold{Z}\,[GL(2,\bold{Z})]\,\widehat{}\ 
$$
where the hat means an obvious  localization
of the group ring. This operator acts on the space
of $s$--densities  $h(dx)^s$. Then:
$$
L(h(dx)^s) = L_sh\, (dx)^s.
$$
and moreover
$$
\int_{[0,1]} f (dx)^{1-s}\cdot L(h(dx)^s) =
\int_{[0,1]} (f(dx)^{1-s}|T) \, h(dx)^s. 
$$
We now take as our heuristic principle the following prescription:

\smallskip
{\it an $L$--invariant $s$--density
is a substitute of a modular form of weight $2-2s$
on the noncommutative modular curve $B$}.

\smallskip

Notice that motivation for adopting this
principle consists of two steps: first, we replace the
action of $\Gamma$ by that of $T$ (in view of (3.2)),
second, we dualize. Roughly speaking, we replace
invariant vectors by invariant functionals. 

\smallskip

The classical example is {\it Gauss 1--density}
$$
\dfrac{1}{\roman{log}\,2}
\dfrac{1}{1+x}
$$
which appeared in the Gauss famous conjecture on the distribution
of continued fractions.

\smallskip

In the following we will briefly describe recent work of
D.~Mayer ([May]), J.~Lewis, and D.~Zagier ([LZ1], [LZ2]), giving a very beautiful
description of the spectrum of values of $s$ for which
an $L$--invariant density exists. The construction starts with
D.~Mayer's discovery of a space on which $L_s$ becomes a 
honest trace class operator.

\medskip

{\bf 3.4. Mayer's operators.} Mayer's space $\bold{V}$ is defined
as the space of holomorphic functions
in $\bold{D}=\{z\in \bold{C}\,|\,|z-1|<\frac{3}{2}\}$
continuous at the boundary. With the supremum norm, it becomes
a complex Banach space.

\medskip

{\bf 3.4.1. Claim.} {\it (i) The formal operator $L_s$ for $\roman{Re}\,s>1/2$
is of trace class (in fact, nuclear of order 0) on the Banach space 
$\bold{V}$.

\smallskip

(ii) It has a meromorphic continuation to the whole complex plane
of $s$, holomorphic except for simple poles at $2s=1,0,-1, \dots$

\smallskip

(iii) The Fredholm determinant $\roman{det}\,(1-L_s^2)$ can be identified with
the Selberg zeta function of $PSL(2,\bold{Z})\setminus H.$}

\smallskip

For a proof, see [May].

\medskip

{\bf 3.4.2. Corollary.} {\it $L$--invariant/antiinvariant $s$--densities which can be obtained by restriction from a density in $\bold{V} |dz|^s$ exist if and only if
$s$ is a zero of the Selberg's zeta $Z(s)$.}

\medskip

All zeroes of $Z(s)$ can be subdivided into following groups:

\smallskip

(i) $s=1$.

\smallskip

(ii) Zeroes on $\roman{Re}\,s=1/2.$

\smallskip

(iii) Zeroes $s=1-k$, $k=2,3,4,\dots $.

\smallskip

(iv) Critical zeroes of Riemann's $\zeta (s)$ divided by two
(hence on $\roman{Re}\,s =1/4$, if one believes the Riemann Hypothesis).

\bigskip

{\bf 3.4.3. Theorem (Lewis--Zagier).} {\it (i) $s=1$ corresponds to the Gauss density $(1+x)^{-1}.$ 

\smallskip

ii)  Zeroes on $\roman{Re}\,s=1/2$ produce all
real analytic $L^2$--invariant $s$--densities on $(0,\infty )$
tending to zero as $x\to +\infty.$
They are automatically holomorphic on $\bold{C} - (-\infty ,0]).$

\smallskip

(iii) Zeroes $s=1-k$, $k=2,3,4,\dots $ produce all polynomial
$L^2$--invariant densities,
which are period functions of modular forms of integral
weight on upper half--plane.

\smallskip

(iv) The  $s$--densities corresponding to
the critical zeroes of Riemann's zeta are analytic continuations of
``half Eisenstein series''
$$
h_s(z)=\sum_{m,n\ge 1}(m(z+1)+n)^{-2s} .
$$
}

A moral of this beautiful story from our perspective is this:
the nontrivial zeroes of $Z(s)$ furnish a spectrum of 
complex fractal dimensions; are there spaces behind them?
 
\newpage

\centerline{\bf \S 4. Fractional dimensions}

\smallskip

\centerline{\bf in homological algebra}

\medskip

{\bf 4.1. Introduction.}  A ``dimension'', or ``weight'',
in homological algebra is simply
a super/subscript of the relevant (co)homology group.
Cohomology groups are invariants of a complex, 
considered as an object of a derived/triangulated category. 
Terms of a complex are routinely graded by integers (at least, up to a shift),
and the differential is of degree $\pm 1$.
Are there situations where we get a fractional numbering?

\smallskip

The answer is positive. Such situations arise in the following way.

\smallskip

A derived category $D(\Cal{C})$ of an abelian category $\Cal{C}$
may have other abelian subcategories $\Cal{C}^{\prime}$
satisfying certain compatibility conditions with the triangulated structure and called ``hearts'' of the respective $t$--structures (see [BeBD]).

\smallskip

For any heart $\Cal{C}^{\prime}$, there is a natural cohomology
functor $H^{\Cal{C}^{\prime}}:\,D(\Cal{C})\to \Cal{C}^{\prime}$.

\smallskip

Objects of such hearts can be considered as ``perverse 
modifications'' of the initial objects of $C$ represented by certain
complexes of objects of $C$. In this way, perverse
sheaves were initially defined via perversity functions
by R.~MacPherson. He has also invented a construction which translates
the algebraic notion of perversity function on triangulated spaces into
a geometric notion of perverse triangulation and revives
the old Euclid's intuition in the context of refined perversity: cf. [Vy].
In this context the dimension remains integral.

\medskip

However, some very common  derived categories $D(\Cal{C})$, for example 
coherent sheaves on an elliptic curve, have {\it families
of hearts $\Cal{C}^{\theta}$ indexed by real numbers $\theta$}.

\smallskip

This ``flow of charges'' was first discovered in the context of
Mirror Symmetry (cf [Dou]).  For a mathematical treatment due to T. Bridgeland,
see [Br]; one can find there some physics comments as well.

\smallskip

The values of the respective homology functor $H^{\theta}:\,D(\Cal{C})\to \Cal{C}^{\theta}$ then
naturally can be thought as having ``dimension/weight $\theta$''.

\smallskip

Below I will describe a particular situation where such groups
appear in the framework of noncommutative tori and
Real Multiplication program of [Ma3] developed in [Po2].

\medskip

{\bf 4.2. $\bold{C}\bold{R}$--lattices.} In the Real
Multiplication program, I suggested to consider 
{\it pseudolattices} which are groups $\bold{Z}^2$
embedded in $\bold{R}$ as ``period lattices of
noncommutative tori'', in the same way as discrete
subgroups $\bold{Z}^2\subset\bold{C}$ are period lattices
of elliptic curves. Here I will start with introducing
the category of {\it $\bold{C}\bold{R}$}--lattices combining properties
of lattices and pseudolattices.

\smallskip

Objects of this category are maps $(j:\,P \to V,s)$, embeddings of
$P\cong \bold{Z}^3$ into 1--dim $\bold{C}$--space $V$,
such that the closure of $j(P)$ is an infinite union of translations
of a real line; $s$ denotes a choice of its orientation.  

\smallskip

(Weak) morphisms are commutative diagrams

$$
\CD 
P^{\prime} @>j^{\prime}>> V^{\prime}\\
@V\varphi VV  @VV\psi V  \\
P @>>j> V  \\
\endCD
$$
where $\psi$ is a linear map.
Strong morphisms should conserve the orientations.

\smallskip

Each  $\bold{C}\bold{R}$--lattice is isomorphic to
one of the form 
$$
P_{\theta ,\tau}:\ \bold{Z}\oplus
\bold{Z}\theta \oplus \bold{Z}\tau \subset \bold{C},\ \theta\in\bold{R},\
\roman{Im}\,\tau \ne 0.
$$
We have $P_{\theta ,\tau}\cong P_{\theta^{\prime},\tau^{\prime}}
\Longleftrightarrow \theta^{\prime}=\dfrac{a\theta +b}{c\theta +d}$
for some
$
\left(\matrix 
a&b\\c&d
\endmatrix\right) \in GL(2,\bold{Z}),
$
and $\tau^{\prime}=\dfrac{\tau +e\theta +f}{c\theta +d}$
for some $e,f\in\bold{Z}$. 

\smallskip

We will discuss the following problem:

\smallskip

{\it Define and study noncommutative
spaces representing quotient spaces $\bold{C}/(\bold{Z}\oplus
\bold{Z}\theta \oplus \bold{Z}\tau )$.}

\medskip

{\bf 4.2.1. Approach I: via noncommutative tori.}
(i) Produce
the quotient $\bold{C}/(\bold{Z}\oplus\bold{Z}\theta )$
interpreting it as a $C^{\infty}$ noncommutative
space, quantum torus $T_{\theta}$,
with the function ring $A_{\theta}$ generated by the unitaries
$$
U_1, U_2,\ U_1U_2=e^{2\pi i\theta}U_2U_1.
\eqno(4.1)
$$

(ii) Introduce ``the complex structure'' on $A_{\theta}$ as
a noncommutative $\overline{\partial}$--operator $\delta_{\tau}$:
$$
\delta_{\tau}U_1:=2\pi i\,\tau U_1,\ \delta_{\tau}U_2:=2\pi i\, U_2 .
\eqno(4.2)
$$
Denote the resulting space $T_{\theta,\tau}$.

\medskip

{\bf 4.2.2. Approach II: via elliptic curves.}
(i) Produce
the quotient $\bold{C}/(\bold{Z}\oplus\bold{Z}\tau )$
interpreting it as the elliptic curve $S_{\tau}$, endowed with a real
point 
$$
x_{\theta}:=\theta\,\roman{mod}\,(\bold{Z}\oplus
\bold{Z}\tau ).
$$

(ii) Form a ``noncommutative quotient space'' $S_{\tau,\theta}:=S_{\tau}/(x_{\theta})$
interpreting it via a crossed product construction
in the language of noncommutative spectra of twisted
coordinate rings, or $q$--deformations of elliptic
functions, etc. 

\medskip

{\bf 4.2.3. Problem.} In what sense the two constructions produce
one and the same noncommutative space $T_{\theta,\tau}$``='' 
$S_{\tau}/(x_{\theta})$?

\medskip

This particular situation is an example of a deeper problem, which more
or less explicitly arises at all turns in noncommutative geometry:
in our disposal there is {\it no good, or even working,} definition
of morphisms of noncommutative spaces. Worse,
we do not quite know what are {\it isomorphisms} between noncommutative spaces.
Even more concretely, assume that our noncommutative space
is a quotient of some common space $M$ with respect to
an action of two commuting groups, say, $G$ and $H$,
such that $G\setminus M$ exists as a honest commutative space
upon which $H$ acts in a ``bad'' way, and similarly
$M/H$ exists as a honest commutative space upon which $G$ 
acts in a ``bad'' way. We can interpret the noncommutative
space $G\setminus M/H$ via a crossed product construction
applied either to $(G\setminus M)/H$ or to $G\setminus (M/H)$.
However, these two constructions generally produce quite
different noncommutative rings.

\smallskip

According to the general philosophy, these two rings then are expected
to be Morita equivalent in an appropriate sense, that is,
have equivalent categories of representations.
For an example where such a statement is a theorem, see
Rieffel's paper [Rie].

\smallskip

Polishchuk's answer to the problem 4.2.3 is in this spirit,
but more sophisticated. We will succinctly state it right away,
and then produce some explanations in the subsections 4.3--4.6:

\medskip

{\bf 4.2.4. Claim (A. Polishchuk).} $T_{\theta,\tau}$ and 
$S_{\tau}/(x_{\theta})$ have canonically equivalent
categories of coherent sheaves which are defined as follows:

\smallskip

On $T_{\theta,\tau}$: the category of vector bundles, that is 
projective modules
over $T_{\theta}$ endowed with a $\delta_{\tau}$--compatible
complex structure.

\smallskip

On $S_{\tau ,\theta}$: the heart $\Cal{C}^{-\theta^{-1}}$ of the $t$--structure of the derived
category of coherent sheaves on the elliptic curve $\bold{C}/
(\bold{Z}+\bold{Z}\tau )$ associated with the slope $-\theta^{-1}$.
 
\medskip

{\bf 4.3. Holomorphic structures on modules and bimodules.}
Let $A_{\theta}$ be the function ring of a $C^{\infty}$ noncommutative torus
(4.1),
endowed with a complex structure (4.2). A right $A_{\theta}$--module $E$
can be geometrically interpreted as vector bundle on this torus.

\smallskip

A holomorphic structure,
compatible with $\delta_{\tau}$, on a right $A_{\theta}$--module $E$ is defined as a map $\overline{\nabla}:\,E\to E$ satisfying
$$
\overline{\nabla}(ea)=\overline{\nabla}(e) a+e\delta_{\tau}(a).
$$
Holomorphic maps between modules are maps commuting with $\overline{\nabla}$.
Similar definitions can be stated for left modules
(see [Co1], [PoS], [Po1], [Po2]).

\smallskip

We imagine modules endowed with such a holomorphic structure as (right)
vector bundles on the respective holomorphic torus.
Their cohomology groups can be defined \`a la Dolbeault: $H^0(E):= \roman{Ker}\,\overline{\nabla}$,
$H^1 (E):= \roman{Coker}\,\overline{\nabla} $.

\smallskip

Similarly, a projective $A_{\theta}$--$A_{\theta^{\prime}}$--bimodule $E$
can be imagined as a sheaf on the product of two smooth tori.
If both tori are endowed with holomorphic structures,
and $E$ is endowed with an operator $\overline{\nabla}$ compatible
with both of them, this sheaf descends to the respective holomorphic
noncommutative space. It can be considered now as
{\it a Morita morphism} $T_{\theta,\tau}\to T_{\theta^{\prime},\tau^{\prime}}$,
that is, the functor $E\otimes * $ from the category
of right vector bundles on $T_{\theta,\tau}$ to the one on 
$T_{\theta^{\prime},\tau^{\prime}}$. This  agrees very well with 
motivic philosophy in commutative algebraic geometry
where morphisms between, say, complete smooth varieties are 
correspondences.

\smallskip

For two--dimensional noncommutative tori, there are good classification
results for these objects.

\medskip

{\bf 4.3.1. $C^{\infty}$--classification of projective $A_{\theta}$--modules.}
Fix $(n,m)\in\bold{Z}^2$. Define the
right $A_{\theta}$--module $E_{n,m}(\theta )$ as 
$A_{\theta}^{|n|}$ for $m=0$ and as the Schwartz space $S(\bold{R}\times \bold{Z}/m\bold{Z})$ with right action
$$
fU_1 (x,\alpha ):= f(x-\frac{n+m\theta}{m},\alpha -1),\
fU_2 (x,\alpha )=\roman{exp}\,(2\pi i(x-\frac{\alpha n}{m})) f(x,\alpha ),
$$
for $m\ne 0.$ These modules are projective.

\smallskip

{\bf 4.3.2. Claim.} {\it  Any finitely generated projective right $A_{\theta}$--module is isomorphic
to $E_{n,m}(\theta )$ with $n+m\theta \ne 0$.}

\smallskip

We define the degree, rank, and slope of $E_{n,m}(\theta )$
by the formulas
$$
\roman{deg}\,E_{n,m}(\theta ):= m,\ \roman{rk}\,E_{n,m}(\theta ):=n+m\theta ,\
\mu\,(E_{n,m}(\theta )):=\frac{m}{n+m\theta}.
$$
Notice that rank is generally fractional: according to the
von Neumann -- Murray philosophy, it is the normalized
trace of a projection.

\smallskip

$E_{n,m}(\theta )$ is isomorphic to $E_{-n,-m}(\theta )$.
It is sometimes convenient to introduce a $\bold{Z}_2$--grading declaring 
$E_{n,m}(\theta )$
even for $\roman{deg}\,E_{n,m}(\theta )>0$ and odd for $<0$.

\smallskip

{\it Basic modules}  are defined as $E_{n,m}(\theta )$ for $(n,m)=1$.
Generally, $E_{nd,md}(\theta ) \cong E_{n,m}(\theta )^d$.

\smallskip

It is convenient to do bookkeeping using
matrices instead of vectors. For $g\in SL(2,\bold{Z})$, write
$E_g(\theta ):=E_{d,c}(\theta )$ where $(c,d)$ is the lower row of $g$.
Accordingly, set
$$
\roman{deg}\, g:= c,\ \roman{rk}\,(g,\theta ):=c\theta +d,
$$
and we have
$$
\roman{rk}\,(g_1g_2,\theta )=\roman{rk}\,(g_1,g_2\theta )\roman{rk}\,(g_2,\theta ),\
g\theta:=\frac{a\theta +b}{c\theta +d}.
$$ 

{\bf 4.3.3. Claim.} {\it The endomorphisms of a basic module  $E_g(\theta )$ form 
an algebra isomorphic to $A_{g\theta}$:
$$
V_1f(x,\alpha )=f(x-\frac{1}{c},\alpha - a),\
V_2f(x,\alpha )=\roman{exp}\,(2\pi i\left(\frac{x}{c\theta + d}-\frac{\alpha}{c}\right))f(x,\alpha).
$$
\smallskip
Thus   $E_g(\theta )$ is  a biprojective $A_{g\theta}$--$A_{\theta}$
bimodule. Tensor multiplication by it produces
a Morita isomorphism of the respective quantum tori.}

\medskip

{\bf 4.4. Holomorphic structures.} Now endow $A_{\theta}$ with the
holomorphic structure (4.2). There is an one--parametric series
of compatible holomorphic structures on  $E_{n,m}(\theta ), m\ne 0$.
Namely, let $z\in \bold{C}$. 
For $f$ in the Schwartz space put
$$
\overline{\nabla}_z(f):=\partial_xf+2\pi i(\tau\mu(E)x+z)f\,.
$$
Similarly,  for $E=A_{\theta }$ as right module put
$$
\overline{\nabla}_z(a):=2\pi iza+\delta_{\tau}(a).
$$
In fact, $\overline{\nabla}_z$ up to isomorphism depends only $z\,\roman{mod}\,
(\bold{Z}+\bold{Z}\tau )/\roman{rk}\,E.$

\smallskip

The Leibniz rule for the left action of $A_{g\theta}$ reads
$$
\overline{\nabla}_z(be)=b\overline{\nabla}_z(e)+\frac{1}{\roman{rk}\,E}
\delta_{\tau}(b)e.
$$

\medskip

{\bf 4.4.1. Theorem.} {\it Right holomorphic bundles on $T_{\theta, \tau}$ 
with arbitrary
complex structures compatible with $\delta_{\tau}$
form an abelian category $\Cal{C}_{\theta,\tau}$.

\smallskip

Every object of this category admits a finite filtration whose
quotients are isomorphic to standard bundles
(basic modules with a standard structure).}

\smallskip

Let us now turn to elliptic curves. We start with some 
cohomological preliminaries.

\medskip

{\bf 4.5. Torsion pairs.} Let $\Cal{C}$ be an abelian category.

\smallskip

{\it A torsion pair in $\Cal{C}$} is
a pair of full subcategories $p=(\Cal{C}_1,\Cal{C}_2)$
stable under extensions, with $\roman{Hom}\,(\Cal{C}_1,\Cal{C}_2)=0$,
and such that every $A\in\Cal{C}$ has a unique subobject in 
$\Cal{C}_1$ with quotient in $\Cal{C}_2$.

\smallskip

(ii) {\it The $t$--structure on $D(\Cal{C})$ associated
with this torsion pair} is defined by
$$
D^{p,\le 0}:= \{K\in D(\Cal{C})\,|\, H^{>0}(K)=0, H^0(K)\in \Cal{C}_1\},
$$
$$
D^{p,\ge 1}:= \{K\in D(\Cal{C})\,|\, H^{<0}(K)=0, H^0(K)\in \Cal{C}_2\}.
$$

\smallskip

(iii) {\it The heart} of this $t$--structure
is a new abelian category $\Cal{C}^p:=D^{p,\le 0}\cap D^{p,\ge 0}$.
It is endowed with a torsion pair as well, namely 
$(\Cal{C}_2[1],\Cal{C}_1)$ ({\it the tilting.})

\smallskip

Let now $X$ be a mooth complete algebraic curve.
{\it Slope} of a stable  bundle on X is defined as  $\roman{deg}/\roman{rk}$,
slope of a torsion sheaf is $+\infty$. The category
$\roman{Coh}_I$, by definition, consists of extensions of sheaves with slope in $I\subset \bold{R}.$ 

\smallskip

For any irrational $\theta$, the pair  $(\roman{Coh}_{>\theta},\roman{Coh}_{<\theta})$
is a torsion pair in the category of coherent sheaves on a curve $X$ ($\theta$ irrational). 

\smallskip

Denote by $\roman{Coh}^{\theta}(X)$ the heart of the respective $t$--structure.

\medskip

{\bf 4.5.1. Theorem.} {\it The category $\Cal{C}_{\theta,\tau}$
defined in the Theorem 4.4.1
is equivalent to
$\roman{Coh}^{-\theta^{-1}}(T_{0,\tau})$.}

\medskip

{\bf 4.6. Real multiplication case.} In the real multiplication case,
$\theta$ is a real quadratic irrationality. The crucial observation
characterising such $\theta$ is this:

\smallskip

{\it $\theta$ is a real quadratic
irrationality $\Longleftrightarrow$

\smallskip

$\exists\ g\in SL(2,\bold{Z})$ such that $g\theta =\theta$
$\Longleftrightarrow$

\smallskip

$E_g(\theta)$ is a biprojective $A_{\theta}$--$A_{\theta}$
bimodule, inducing a nontrivial autoequivalence of
the category of $A_{\theta}$--modules.}

\medskip

Polishchuk remarked that  such a bimodule, endowed with
a holomorphic structure, determines
a noncommutative graded ring   
$$
B:= \oplus_{n\ge 0} H^0(E_g(\theta)^{\otimes n}_{A_{\theta}}).
$$
This ring is a particular case of a more general
categorical construction. Let $\Cal{C}$ be an additive category,
$F:\,\Cal{C} \to \Cal{C}$ an additive functor, $O$ an object.

\smallskip

Then we have a  graded ring with twisted multiplication:
$$
A_{F,O}:=\oplus_{n\ge 0} \roman{Hom}_{\Cal{C}}(O,F^n(O)).
$$

\smallskip

To produce the former $B$, choose $O=A_{\theta}$, $F(\cdot )=\cdot\otimes E_g(\theta )$
in the category of holomorphic bundles.

\medskip

{\bf 4.6.1. Theorem (A.~Polishchuk).} {\it For every real quadratic irrationality $\theta$
and complex structure $\tau$, the heart $\roman{Coh}^{\theta}(T_{0,\tau})$
is equivalent to the Serre category of right coherent sheaves
on the noncommutative projective spectrum of some algebra of the form
$A_{F,O}$ with $F:\,D^b(T_{0,\tau})\to D^b(T_{0,\tau})$
being an autoequivalence.}

\medskip

This remarkable result may be tentatively considered as an approximation
to the problem invoked in [Ma4]: how to find finitely
generated (over $\bold{C}$ and eventually over $\bold{Z}$)
rings naturally associated with $T_{\theta, \tau}$. They are necessary
to do arithmetics.

\smallskip

Since however an arbitary parameter $\tau$ appears in this
construction, the following question remains:
how to choose $\tau$ for arithmetical
applications? 

\smallskip

A natural suggestion is: if $\theta\in \bold{Q}(\sqrt{d})$, choose
$\tau\in \bold{Q}(\sqrt{-d})$. 

\smallskip

For such a choice, Polishchuk's noncommutative rings can be considered
as a more sophisticated version of Kronecker's idea
to merge $\sqrt{d}$ with $\sqrt{-d}$ in order to produce
solutions of Pell's equation in terms of elliptic functions:
cf. A.~Weil [We] for a modern exposition and historical
context.

\bigskip

\centerline{\bf Bibliography}

\medskip

[And] G.~Anderson. {\it Cyclotomy and a covering of the Taniyama group.}
Comp. Math., 57 (1985), 153--217.

\smallskip

[At] M.~Atiyah. {\it Commentary on the article of Yu.~I.~Manin:
``New dimensions in geometry.''} In: Springer Lecture
Notes in Math., 1111, 1985, 103--109.

\smallskip

[BeBD] A.~Beilinson, J.~Bernstein, P.~Deligne. {\it Faisceaux pervers.}
Ast\'erisque, 100 (1982), 5--171.

\smallskip

[Br] T.~Bridgeland. {\it Stability conditions on triangulated
categories.} Preprint math.AG/0212237

\smallskip

[Co1] A.~Connes. {\it $C^*$ alg\`ebres et g\'eom\'etrie
differentielle.} C. R. Acad. Sci. Paris, S\'er. A--B, 290 (1980),
A599--A604.

\smallskip

[Co2] A.~Connes. {\it Noncommutative Geometry.} Academic Press, 1994.

\smallskip

[Co3] A.~Connes, {\it Geometry from the spectral point of view}.
Lett. Math. Phys. Vol. 34 (1995), 203--238.

\smallskip

[Co4] A.~Connes, {\it Gravity coupled with matter and the
foundation of non-commuta\-tive geometry}.  Comm. Math. Phys.  182
(1996),  no. 1, 155--176.

\smallskip  

[Co5] A.~Connes, {\it Trace formula in noncommutative
geometry and the zeros of the Riemann zeta function}.  Selecta
Math. (N.S.)  5  (1999),  no. 1, 29--106. 

\smallskip 

[Co6] A.~Connes, {\it Cyclic cohomology, quantum group
symmetries and the local index formula for $SU_q(2)$},  
J. Inst. Math. Jussieu  3  (2004),  no. 1, 17--68.

[CoKr] A.~Connes, D.~Kreimer. {\it Renormalization in
quantum field theory and the Riemann--Hilbert problem.
I. The Hopf algebra structure of graphs and the main theorem.}
Comm. Math. Phys. 210:1 (2000), 249--273.

\smallskip

[CoMar1] A.~Connes, M.~Marcolli. {\it Quantum statistical mechanics
of $\bold{Q}$--lattices. (From Physics to Number Theory
via Noncommutative Geometry, Part I).} Preprint math.NT/0404128

\smallskip

[CoMar2] A.~Connes, M.~Marcolli. {\it Renormalization and
motivic Galois theory}. Int. Math. Res. Notices ,
N.76 (2004), 4073--4092.

\smallskip

[CoMos] A.~Connes, H.~Moscovici. {\it The local index formula
in noncommutative geometry}. GAFA Vol.5 (1995), N.2, 174--243.

\smallskip

[ConsMar] C.~Consani, M.~Marcolli, {\it Noncommutative geometry,
dynamics, and $\infty$-adic Arakelov geometry}.  Selecta Math. (N.S.)
10  (2004),  no. 2, 167--251.

\smallskip

[DaLSSV] L.~D\c{a}browski, G.~Landi, A.~Sitarz, W.~van Suijlekom,
J.C.~Varilly, {\it The Dirac operator on $SU_q(2)$}. Preprint math.QA/0411609.

\smallskip

[De1] C.~Deninger. {\it Local $L$--factors of motives and regularized determinants}. Inv. Math. 107 (1992), 135--150.

\smallskip

[De2] C.~Deninger. {\it Motivic $L$--functions and regularized determinants.}
Proc. Symp. Pure Math., 55:1 (1994), 707--743

\smallskip

[De3] C.~Deninger. {\it Some analogies between number theory and dynamical systems on foliated spaces.} Doc. Math. J. DMV. Extra volume
ICM I (1998), 23--46.

\smallskip

[De4] C.~Deninger. {\it A note on arithmetic topology and dynamical systems.}
Contemp. Math. 300, AMS, Providence RI (2002), 99--114.

\smallskip

[Dou] M.~Douglas. {\it Dirichlet branes, homological mirror symmetry, and stability.} In: Proceedings of the ICM 2002, Higher Education Press,
Beijing 2002, vol. III, 395--408.   Preprint math.AG/0207021

\smallskip

[Ha] M.~Haran. {\it The mysteries of the real prime.}
Clarendon Press, Oxford, 2001.

\smallskip

[He] Heath T. L. {\it The thirteen books of Euclid's
Elements. Translation, Introduction and Commentary.} Cambridge UP, 1908.

\smallskip

[KalW] W.~Kalau, M.~Walse, {\it Gravity, non-commutative
geometry and the Wodzicki residue}, J. Geom. Phys. Vol.16 (1995),
327--344.

\smallskip

[KapSm] M.~Kapranov, A.~Smirnov. {\it Cohomology determinants
and reciprocity laws: number field case.} Unpublished.

\smallskip

[Ka] N. Katz. {\it $p$--adic properties of modular schemes and modular forms.}
Springer LNM, 350 (1973)

\smallskip

[Kr] D.~Kreimer. {\it On the Hopf algebra structure
of perturbative Quantum Field Theory.} Adv. Theor. Math. Phys.,
2:2 (1998), 303--334.

\smallskip

[LaPo] M.~Lapidus, C.~Pomerance. {\it Fonction z\^eta de Riemann
et conjecture de Wyl--Berry pour les tambours fractals.}
C. R. Ac. Sci. Paris, S\'er. I Math., 310 (1990), 343--348.

\smallskip

[LavF1] M.~Lapidus, M.~van Frankenhuysen. {\it Complex dimensions
of fractal strings and oscillatory phenomena in fractal
geometry and arithmetic.}  In: Spectral Problems
in Geometry and Arithmetic (ed. by T.~Branson), Contemp.
Math., vol. 237, AMS, Providence R.I., 1999, 87--105.

\smallskip

[LavF2] M.~Lapidus, M.~van Frankenhuysen. {\it Fractal Geometry 
and number theory. Complex dimensions of fractal strings and zeros of zeta functions.} Birkh\"auser, 1999.

\smallskip

[LZ1] J. Lewis, D. Zagier. {\it Period functions for Maass wave forms.}
Ann. of Math., 153 (2001), 191--258.

\smallskip

[LZ2] J. Lewis, D. Zagier. {\it Period functions and the Selberg zeta function
for the modular group.} In: The mathematical beauty
of physics, Adv. Series in Math. Physics, 24, World Sci. Publ.,
River Edge, NJ (1997), 83--97.

\smallskip

[Mand] B. Mandelbrot. {\it The fractal geometry of nature.} Freeman \& Co, NY, 1983.

\smallskip

[Ma] D.~Yu.~Manin. {\it Personal communication.}

\smallskip

[Ma1] Yu.~Manin. {\it New dimensions in geometry.} Russian:
Uspekhi Mat. Nauk, 39:6 (1984), 47--73. English:
Russian Math. Surveys, 39:6 (1984), 51--83,
and  Springer Lecture
Notes in Math., 1111, 1985, 59-101.

\smallskip

[Ma2] Yu.~Manin. {\it Three--dimensional hyperbolic geometry as $\infty$--adic
Arakelov geometry.}  Inv. Math., 104 (1991), 223--244.

\smallskip

[Ma3] Yu.~Manin. {\it Lectures on zeta functions and motives (according to Deninger and Kurokawa).}
In: Columbia University Number Theory Seminar,  Ast\'erisque 228
(1995), 121--164.

\smallskip

[Ma4] Yu.~Manin. {\it Von Zahlen und Figuren.} Preprint math.AG/0201005

\smallskip

[Ma5] Yu.~Manin. {\it Real multiplication and noncommutative
geometry.} In: The legacy of Niels Henrik Abel,
ed. by O.~A.~Laudal and R.~Piene, Springer Verlag,
Berlin 2004, 685--727. Preprint math.AG/0202109

\smallskip

[MaMar] Yu.~Manin, M.~Marcolli. {\it Continued fractions, modular symbols, and non-commutative geometry.} Selecta math., new ser. 8 (2002),
475--521. Preprint math.NT/0102006

\smallskip

[May] D.~Mayer. {\it Continued fractions and related
transformations.} In: Ergodic Theory, Symbolic Dynamics
and Hyperbolic Spaces, Eds. T.~Bedford et al., Oxford
University Press, Oxford 1991, pp. 175--222.

\smallskip

[Maz] B.~Mazur. {\it Note on \'etale cohomology of number fields.}
Ann. Sci. ENS 6 (1973), 521--552.

\smallskip

[Mor] M.~Morishita. {\it On certain analogies between
knots and primes.} Journ. f\"ur die reine u. angew. Math.,
550 (2002), 141--167.

\smallskip

[PoS] A.~Polishchuk, A.~Schwarz. {\it Categories of holomorphic
bundles on noncommutative two--tori.} Comm. Math. Phys.
236 (2003), 135--159. Preprint math.QA/0211262

\smallskip

[Po1] A.~Polishchuk. {\it Noncommutative two--tori with real
multiplication as noncommutative projective varieties.}
Preprint math.AG/0212306

\smallskip

[Po2]  A.~Polishchuk. {\it Classification of holomorphic vector bundles on noncommutative two--tori.} Preprint math.QA/0308136

\smallskip 

[Rie] M.~Rieffel. {\it Von Neumann algebras associated
with pairs of lattices in Lie groups.} Math. Ann. 257 (1981),
403--418.

\smallskip

[Se] J.--P. Serre. {\it Formes modulaires et fonctions z\^eta $p$--adiques.}
Springer LNM, 350 (1973), 191--268

\smallskip

[Sm1] A.~Smirnov. {\it Hurwitz inequalities for number fields.}
St. Petersburg Math. J., 4 (1993), 357--375.

\smallskip

[Sm2] A.~Smirnov. {\it Letters to Yu. Manin of Sept. 29 and Nov. 29, 2003.}

\smallskip

[Sou] C.~Soul\'e. {\it Les vari\'et\'es sur le corps \`a un el\'ement.}
Moscow Math. Journ., 4:1 (2004), 217--244.

\smallskip

[Ta] J.~Tate. {\it Duality theorems in Galois cohomology over number fields.}
Proc. ICM, Stockholm 1962, 288--295.

\smallskip

[Ti] J.~Tits. {\it Sur les analogues alg\'ebriques des
groupes semi--simples complexes.} Colloque d'Alg\`ebre
sup., Bruxelles 1956. Louvain, 1957, 261--289.

\smallskip

[Vy] M.~Vybornov. {\it Constructible sheaves on simplicial
complexes and Koszul duality.} Math. Res. Letters,
5 (1998), 675--683.

\smallskip

[We] A.~Weil. {\it Elliptic functions according to Eisenstein
and Kronecker.} Springer Verlag, Berlin, 1976 and 1999.

\enddocument